\documentclass[12pt]{article}
\title{Multi-material topology optimization of an electric machine considering demagnetization}

\usepackage{a4,epsfig}
\usepackage{amsfonts,amsmath}
\usepackage{subcaption}
\usepackage{algorithm}
\usepackage{algpseudocode}
\usepackage{pgfplotstable}
\usepackage{hyperref}

\newtheorem{remark}{Remark}[section]

\numberwithin{equation}{section} 

\begin{document}

\setcounter{page}{1}

\title{Multi-material topology optimization of an electric machine considering demagnetization}
\author{Peter~Gangl$^1$, Nepomuk~Krenn$^1$}
\date{$^1$Johann Radon Institute for Computational and
  Applied Mathematics, \\
Altenberger Stra{\ss}e 69, 4040 Linz, Austria}

\maketitle

\begin{abstract}
    We consider the topology optimization problem of a 2d permanent magnet synchronous	machine in magnetostatic operation with demagnetization. This amounts to a PDE-constrained multi-material design optimization problem with an additional pointwise state constraint. Using a generic framework we can incorporate this additional constraint and compute the corresponding topological derivative. We present and discuss optimization results obtained by a multi-material level set algorithm.
\end{abstract}

\section{Introduction}
\begin{figure}
	\centering
	\includegraphics[trim=0.25cm 0cm 2.7cm 0cm,clip,height=5.1cm]{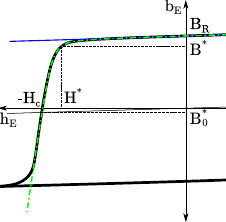}\qquad
	\includegraphics[height=5.4cm]{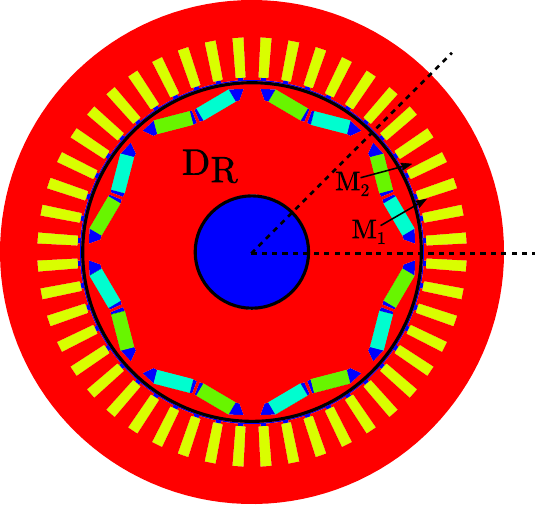}
	\caption{Demagnetization curve in easy direction of permanent magnet with linear (blue, solid) and non-linear (green, dot-dashed) approximation (left). Reference machine consisting of iron in red, air in blue, magnet $M_1$ in light blue and magnet $M_2$ in green (right).}\label{krenn:fig_bhcurve}
\end{figure}

Topology optimization of electric machines is of high interest in industry and research, first introduced in \cite{lowther}. One crucial criterion in the design process of an electric machine is its robustness regarding demagnetization which is an irreversible damage of the permanent magnets. This happens if strong demagnetizing fields occur in direction of magnetization (easy axis) with $h_E\le H^*$, see Fig. \ref{krenn:fig_bhcurve} (left). To avoid this we demand the magnetic flux to fulfill $b_E\ge B^*$. In \cite{theodor} the authors proposed in a first step to model the material law of the permanent magnet by a non-linear relation (green, dot-dashed) instead of using the common ideal linear curve (blue, solid), which leads to designs robust to demagnetization in nominal operation. In \cite{demag} the authors compared different magnet shapes regarding demagnetization. We propose a new way on how to consider demagnetization for an additional damaging operation point by adding a state constraint $b_E^D\ge B^*$ to the topology optimization problem denoting "damage" with the superscript "D".\newline 
We are interested in using the topological derivative for the optimization, a pointwise sensitivity of a cost functional to material changes which was first applied to electric machines in \cite{gangl_diss}. In the multi-material context the topological derivative is vector-valued since various material changes are possible. We use it to update a vector-valued level set representation of the design \cite{gangl_multi},\cite{amstutz_levelset}. To translate the state constraint into a penalty functional which fits into the framework \cite{gangl_auto} to derive the topological derivative formula we apply the approach from \cite{amstutz_state}. \newline The rest of this paper is structured as follows: In Section \ref{krenn:sec_1} we will introduce the machine of interest and the material laws including both a linear and a non-linear one for the permanent magnets. Then we will introduce the new state constraint and provide the topological derivative formula in Section \ref{krenn:sec_2}. Finally in Section \ref{krenn:sec_3} we compare the demagnetization robustness of three designs obtained using a linear magnet model, a non-linear one and the additional state constraint. 
\section{Problem formulation}\label{krenn:sec_1}
We consider the optimization of the rotor $D_R$ of a permanent magnet synchronous machine $D$ with 48 stator slots and 4 pole pairs operating at $15.000\text{rpm}$ displayed in Fig. \ref{krenn:fig_bhcurve} (right). Due to symmetry it is sufficient to optimize only one eighth of the machine with two magnets imposing periodicity conditions on the design. Since we consider only the motor operation and neglect the generator one we will obtain non-symmetric designs. We have ferromagnetic iron in $\Omega_I$, permanent magnets with magnetization direction $e_{M_1}:=\begin{pmatrix}
	\cos\varphi_1\\\sin\varphi_1
\end{pmatrix}, \varphi_1=30^\circ$ in $\Omega_{M_1}$, and $e_{M_2}:=\begin{pmatrix}
	\cos\varphi_2\\\sin\varphi_2
\end{pmatrix}, \varphi_2=15^\circ$ in $\Omega_{M_2}$ and air in $\Omega_A$. We use the notation $\Omega=(\Omega_I,\Omega_{M_1},\Omega_{M_2},\Omega_A)$ to describe the material distribution within the rotor $D_R=\dot{\bigcup}_{i\in\{I,M_1,M_2,A\}}\Omega_i$.
\subsection{Torque maximization} 
The PDE-constrained optimization problem reads
\begin{align}
	\min_{\Omega}\mathcal{T}(\Omega)&=-T(a_\Omega) \text{ subject to } \label{krenn:PF_Tork}\\
	\begin{split}\label{krenn:PF_magneto}
		\widetilde{\text{curl}}h_\Omega(\text{curl} a_\Omega)&=j\text{ in }D\\
		a_\Omega&= 0\text{ on }\partial D
	\end{split}
\end{align}
where the constitutive law $h_\Omega(b):=h_I(b)\chi_{\Omega_I}+h_{M_1}(b)\chi_{\Omega_{M_1}}+h_{M_2}(b)\chi_{\Omega_{M_2}}+h_A(b)\chi_{\Omega_A}$ is defined piecewise. The objective $T(a_\Omega)$ is the torque for one rotor position where $a_\Omega$ denotes the $z$-component of the magnetic vector potential solving the equations of magnetostatics (\ref{krenn:PF_magneto}) with a three-phase sinusoidal excitation with effective value $\hat{j}=1512.5\text{A}$ and an angle $\theta_0=6^\circ$ between the initial electrical and mechanical configuration. Note that the operator $\text{curl} v=\begin{pmatrix}
	\partial_yv\\-\partial_xv
\end{pmatrix}$ is the scalar-to-vector curl, whereas $\widetilde{\text{curl}}v=-\partial_y v_1+\partial_x v_2$ is the vector-to-scalar curl.
\begin{remark}
	Although we stated the optimization problem (\ref{krenn:PF_Tork})-(\ref{krenn:PF_magneto}) for a fixed rotor position we will consider the average torque for $N=11$ equidistant rotor positions. We use harmonic mortaring \cite{egger_mortaring} at the airgap for the efficient coupling of rotor and stator and the consistent torque computation based on the time derivative of the magnetic energy. Since the torque is evaluated in the airgap which is not part of the design domain one can use the topological derivative formula of the torque \cite{gangl_diss} with minor adaptions at the adjoint equation. Nevertheless one could also use other methods for the torque computation.
\end{remark}
\begin{remark}
	One has to add a magnet volume constraint
	\begin{align}\label{krenn:magvol}
		|\Omega_{M_1}\cup\Omega_{M_2}|\le V^*,
	\end{align}
	otherwise the optimal rotor design will consist out of magnets only.
\end{remark}
\subsection{Constitutive laws}
For air we have the relation $h_A(b)=\nu_0b, \nu_0=\frac{10^7}{4\pi}$, for iron we consider an isotropic material law inspired by measurement data
\begin{align*}
	h_I(b)=f_I(|b|)\frac{b}{|b|} \text{ with }f_I(s)=\nu_0s+(200-\nu_0)\frac{2.2s}{\sqrt[12]{2.2^{12}+s^{12}}}.
\end{align*}
For permanent magnets we consider a general material law of the form
\begin{align*}
	h_{M_i}(b)=\begin{pmatrix}
		\cos\varphi_i&-\sin\varphi_i\\\sin\varphi_i&\cos\varphi_i
	\end{pmatrix}\begin{pmatrix}
		f_M(b\cdot e_{M_i})\\
		\nu_mb\cdot e_{M_i}^\perp
	\end{pmatrix}
\end{align*}
for $i=1,2$ with $\nu_m=\frac{\nu_0}{1.086}$ where we choose
\begin{align}\label{krenn:law_linear}
	f_M(s)=\nu_m (s - B_R)
\end{align}
with the remanence flux density $B_R=1.2\text{T}$ for a linear magnet or
\begin{align}\label{krenn:law_nonlinear}
	f_M(s) = \nu_ms+(\frac{\nu_m}{70}-\nu_m)\frac{0.56s}{\sqrt[12]{0.56^{12}+s^{12}}}-H_c
\end{align}
with the coercivity $H_c=4.75\cdot 10^5 \frac{\text{A}}{\text{m}}$ for a non-linear law mimicking the behavior of a N42 magnet at $100^\circ\text{C}$ \cite{arnoldmagnetics}.
\begin{remark}
	The linear material law \eqref{krenn:law_linear} is equivalent to the well-known relation
	\begin{align*}
		h_M(b)=\nu_m(b-B_R).
	\end{align*}
\end{remark}
\section{Demagnetization constraint}\label{krenn:sec_2}
\begin{figure}
	\centering
	\includegraphics[width=7cm]{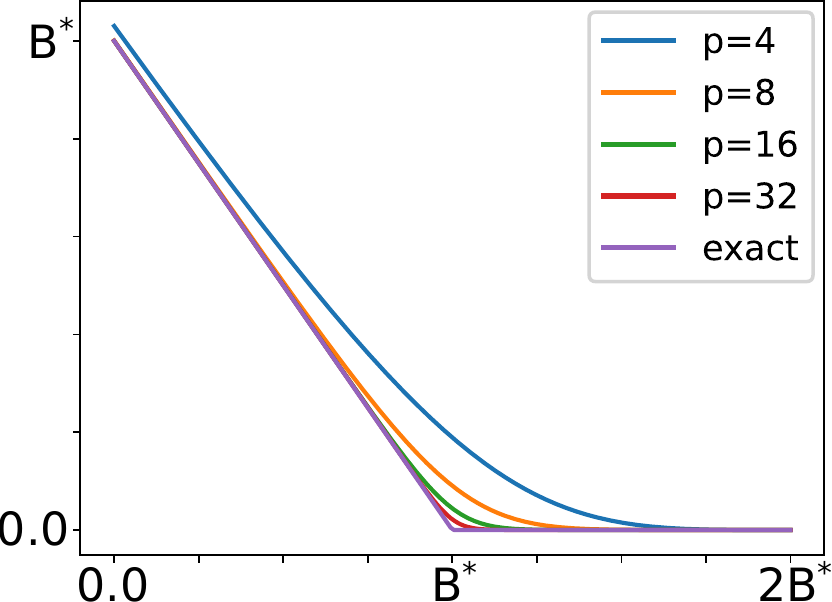}
	\caption{Functions $\Phi$ and $\Phi_{p}$ for $p=4,8,16,32$}
	\label{krenn:fig_phis}
\end{figure}

To avoid demagnetization at a (possibly) different damaging operation point we introduce the constraint
\begin{align}
	\text{curl} a_\Omega^D\cdot e_{M_i}&\ge B^*\text{ a.e. in } \Omega_{M_i}, i=1,2 \text{ where }\label{krenn:PF_constr}\\ \label{krenn:PF_magnetoD}
	\begin{split}
		\widetilde{\text{curl}}h_\Omega(\text{curl} a_\Omega^D)&=j^D\text{ in }D\\
		a_\Omega^D&=0\text{ on } \partial D
	\end{split}
\end{align}
with $B^*=0.56\text{T}$. We transform the state constraint (\ref{krenn:PF_constr}) into a penalty term similarly as in \cite{amstutz_state}:
\begin{align*}
	\text{curl} a_\Omega^D\cdot e_{M_i}\ge B^*\text{ a.e. in }\Omega_{M_i}&\Leftrightarrow B^*-\text{curl} a_\Omega^D\cdot e_{M_i}\le 0\text{ a.e. in }\Omega_{M_i}\\
	\Leftrightarrow\max\left\{B^*-\text{curl} a_\Omega^D\cdot e_{M_i},0\right\}=0&\text{ a.e. in }\Omega_{M_i}\Leftrightarrow\int_{\Omega_{M_i}}\Phi\left(\text{curl} a_\Omega^D\cdot e_{M_i}\right)\text{d}x=0
\end{align*}
for $i=1,2$ with the penalty function $\Phi(s)=\max\{B^*-s,0\}=\max\{2B^*-s,B^*\}-B^*$. Since we aim to compute the topological gradient for this problem we need to approximate $\Phi$ by a smooth function $\Phi_{p}$. Inspired by \cite{amstutz_state} we take
\begin{align}\label{krenn:eq_rho}
	\Phi_{p}(s):=\begin{cases}
		(B^{*p}+(2B^*-s)^p)^\frac{1}{p}-B^*\quad &s<2B^*\\
		0 &s\ge 2B^*
	\end{cases}
\end{align} which regularizes the maximum by the $p$-norm visualized in Fig. \ref{krenn:fig_phis}. One has to choose $p$ carefully since too low values would penalize fields which fulfill the constraint and too high values would result in ill-conditioned behavior of the problem. For this we take $p=16$. 
We get the relaxed constraint
\begin{align}\label{krenn:weak_constr}
	C(a_\Omega^D)=\int_{D_R}j_{\Omega,p}\left(\text{curl} a_\Omega\right)\text{d}x 
\end{align}
with $j_{\Omega,p}(b):=\Phi_p(b\cdot e_{M_1})\chi_{\Omega_{M_1}}+\Phi_p(b\cdot e_{M_2})\chi_{\Omega_{M_2}}$.
With a sufficiently large weight $\gamma$ we get the reduced formulation of the PDE-constrained optimization problem
\begin{align*}\begin{split}
		\min_{\Omega}\mathcal{J}(\Omega)=\mathcal{T}(\Omega)+\gamma \mathcal{C}&(\Omega)+\mathcal{L}(\Omega),\text{ where }
	\end{split}\\
	\mathcal{T}(\Omega)=-T(a_\Omega), a_\Omega \text{ solution of }(\ref{krenn:PF_magneto}),&\ \mathcal{C}(\Omega)=C(a_\Omega^D), a_\Omega^D \text{ solution of }(\ref{krenn:PF_magnetoD})
\end{align*}
and $\mathcal{L}(\Omega)$ is the Augmented Lagrangian \cite{AL} incorporation of the magnet volume constraint \eqref{krenn:magvol}.

\subsection{Topological derivative}
The topological derivative for the torque \cite{gangl_diss} and the volume \cite{amstutz_levelset} are available, we need to derive the formula for the new constraint functional $\mathcal{C}(\Omega)$. With $B:=\text{curl} a^D_\Omega(z)\in\mathbf{R}^2, \hat{B}:=\text{curl} p^D_\Omega(z)\in\mathbf{R}^2$ where the adjoint state $p^D_\Omega\in H^1_0(D)$ is the solution of
\begin{align*}
	\int_D h_\Omega'(\text{curl} a_\Omega^D)\text{curl} p_\Omega^D\cdot\text{curl} v\text{d}x = -\int_{D_R}j'_{\Omega,p}(\text{curl} a_\Omega^D)(\text{curl} v)\text{d}x
\end{align*}
for all $v\in H^1_0(D)$ we can apply the framework \cite{gangl_auto} to get the sensitivity in a point $z\in \Omega_0$ subject to a change from material $0$ to material $1$:
\begin{align}\begin{split}\label{krenn:TD}
		\text{d}^{0\rightarrow 1}\mathcal{C}(\Omega)(z)&=\frac{1}{\pi}\int_{\mathbf{R}^2}(h_\omega(B+\text{curl} K)-h_\omega(B)-h'_\omega(B)(\text{curl} K))\cdot \hat{B}\text{d}\xi \\&+\frac{1}{\pi}\int_\omega (h'_1(B)-h'_0(B))(\text{curl} K)\cdot \hat{B}\text{d}\xi+(h_1(B)-h_0(B))\cdot\hat{B}\\
		&+\frac{1}{\pi}\int_{\mathbf{R}^2}\left(j_{\omega,p}(B+\text{curl} K)-j_{\omega,p}(B)-j_{\omega,p}'(B)(\text{curl} K)\right)\text{d}\xi\\
		&+\frac{1}{\pi}\int_\omega(j_{1,p}'(B)-j_{0,p}'(B))(\text{curl} K)\text{d}\xi+j_{1,p}(B)-j_{0,p}(B)
	\end{split}
\end{align}
where $h_\omega(b):=h_1(b)\chi_\omega+h_0(b)\chi_{\mathbf{R}^2\setminus\overline{\omega}}, j_{\omega,p}(s):=j_{1,p}(s)\chi_\omega+j_{0,p}(s)\chi_{\mathbf{R}^2\setminus\overline{\omega}}$ are piecewise defined with $\omega=B_1(0)$ the open unit ball and the directional derivative of the cost function $j_{\Omega,p}'(b)(\tilde{b})=\Phi_p'(b\cdot e_{M_1})\tilde{b}\cdot e_{M_1}\chi_{\Omega_{M_1}}+\Phi_p'(b\cdot e_{M_2})\tilde{b}\cdot e_{M_2}\chi_{\Omega_{M_2}}$.
The function $K\in\mathcal{B}:=\{v\in L_2(\mathbf{R}^2):\text{curl} v\in L^2(\mathbf{R}^2)\}/\mathbf{R}$ in (\ref{krenn:TD}) is the solution of the auxiliary exterior problem for all $v\in\mathcal{B}$:
\begin{align}\label{krenn:exterior}
	\int_{\mathbf{R}^2}\left(h_\omega(B+\text{curl} K)-h_\omega(B)\right)\cdot \text{curl} v\text{d}\xi = - \int_\omega\left(h_1(B)-h_0(B)\right)\cdot\text{curl} v\text{d}\xi
\end{align}
\subsection{Efficient evaluation}
In every optimization step we need to compute the topological derivative (\ref{krenn:TD}) for all points $z\in D_R\setminus\bigcup_{i\in\{I,M_1,M_2,A\}}\partial \Omega_i$ which depends only on the solution of the state and adjoint equation $B=\text{curl} a_\Omega^D(z),\hat{B}=\text{curl} p_\Omega^D(z)$. This is very costly since we have to solve the non-linear exterior problem (\ref{krenn:exterior}) for every $B$. Therefore we pre-compute the formula for selected samples in an offline phase and evaluate their interpolation in the online phase. Since (\ref{krenn:TD}) is affine linear in $\hat{B}$ it is sufficient to take samples for $B\in\mathbf{R}^2$ in a feasible range.
\begin{remark}
	For isotropic materials one can show rotational invariance of the topological derivative and it is enough to sample for $|B|\in\mathbf{R}_+$, see \cite{gangl_diss}. However, as soon as one considers permanent magnets this is not possible any more.
\end{remark}
\section{Optimization results}\label{krenn:sec_3}
We measure the partial demagnetization of a design by 
\begin{align*}
	\mathcal{D}(\Omega)=\sum_{i=1}^{2}\frac{1}{|\Omega_{M_i}|}\int_{\Omega_{M_i}}\max\left\{\frac{B^*-\text{curl} a_\Omega^D\cdot e_{M_i}}{B^*-B^*_0},0\right\}\text{d}x,
\end{align*}
where $B_0^*=-0.66\text{T}$ indicates the fully demagnetized state. If there is no additional damaging operation point we simply set $a_\Omega^D=a_\Omega$ to measure the demagnetization in nominal operation. For the damaging operation point we choose a high source current $j^D=1.5j$, however considering also different load angles or short circuit operation would be of interest. 
\begin{figure}
	\centering
	\begin{tabular}{ccc}
		linear magnets&non-linear magnets&non-linear magnets with constraint\\ \hline&&\\
		\includegraphics[trim=14.7cm 5.3cm 14.7cm 5.3cm,clip,width=0.32\textwidth]{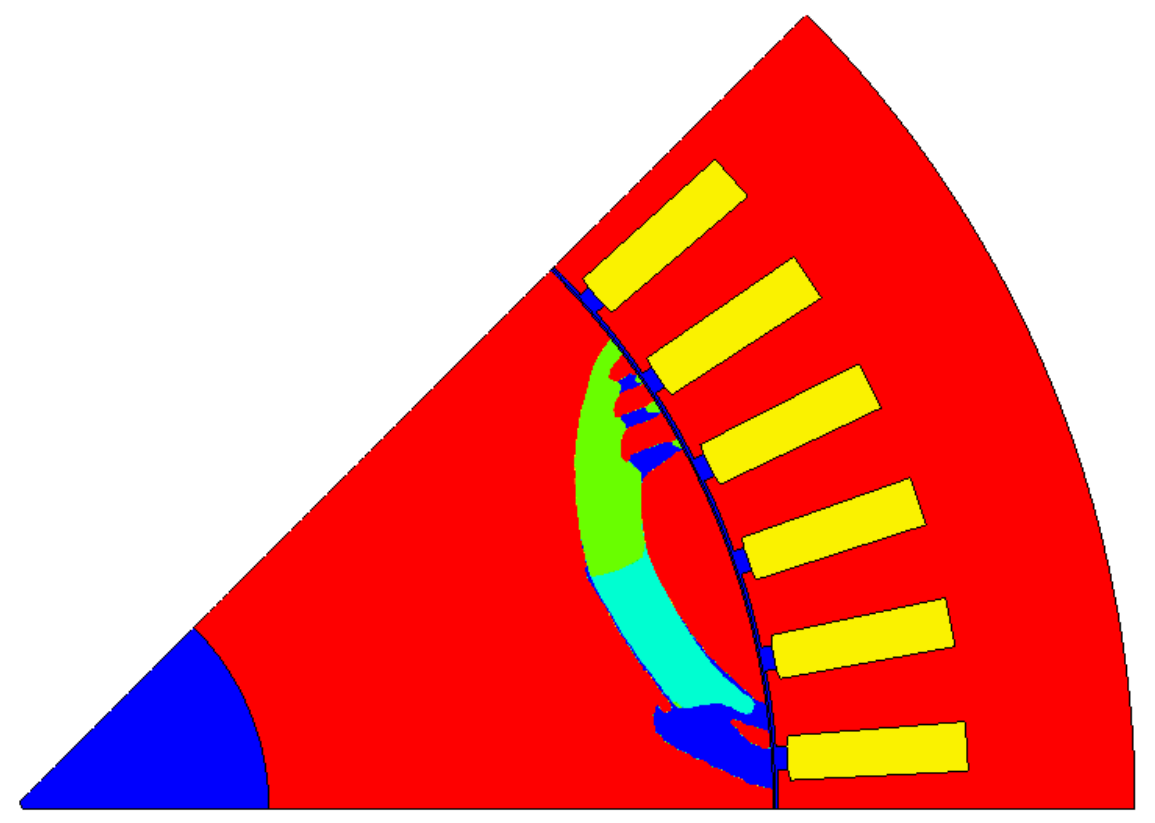} 
		&\includegraphics[trim=14.7cm 5.3cm 14.7cm 5.3cm,clip,width=0.32\textwidth]{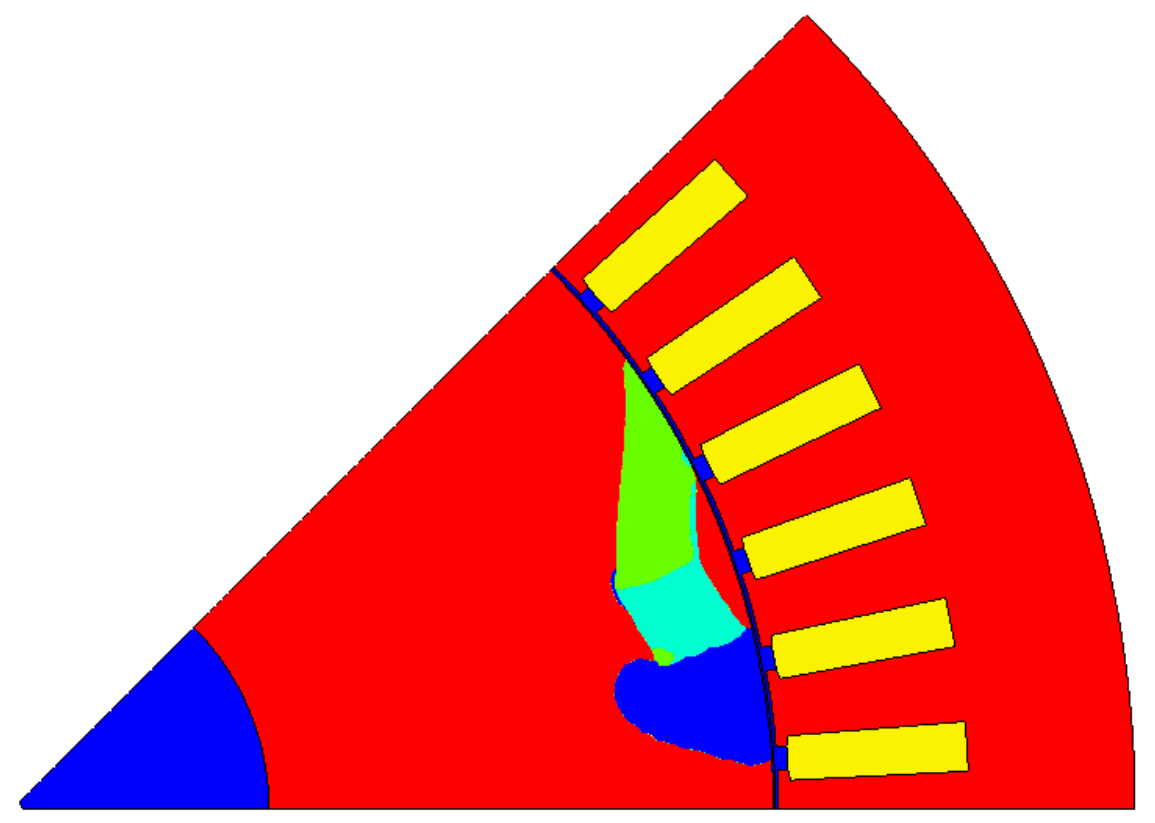} 
		&\includegraphics[trim=14.7cm 5.3cm 14.7cm 5.3cm,clip,width=0.32\textwidth]{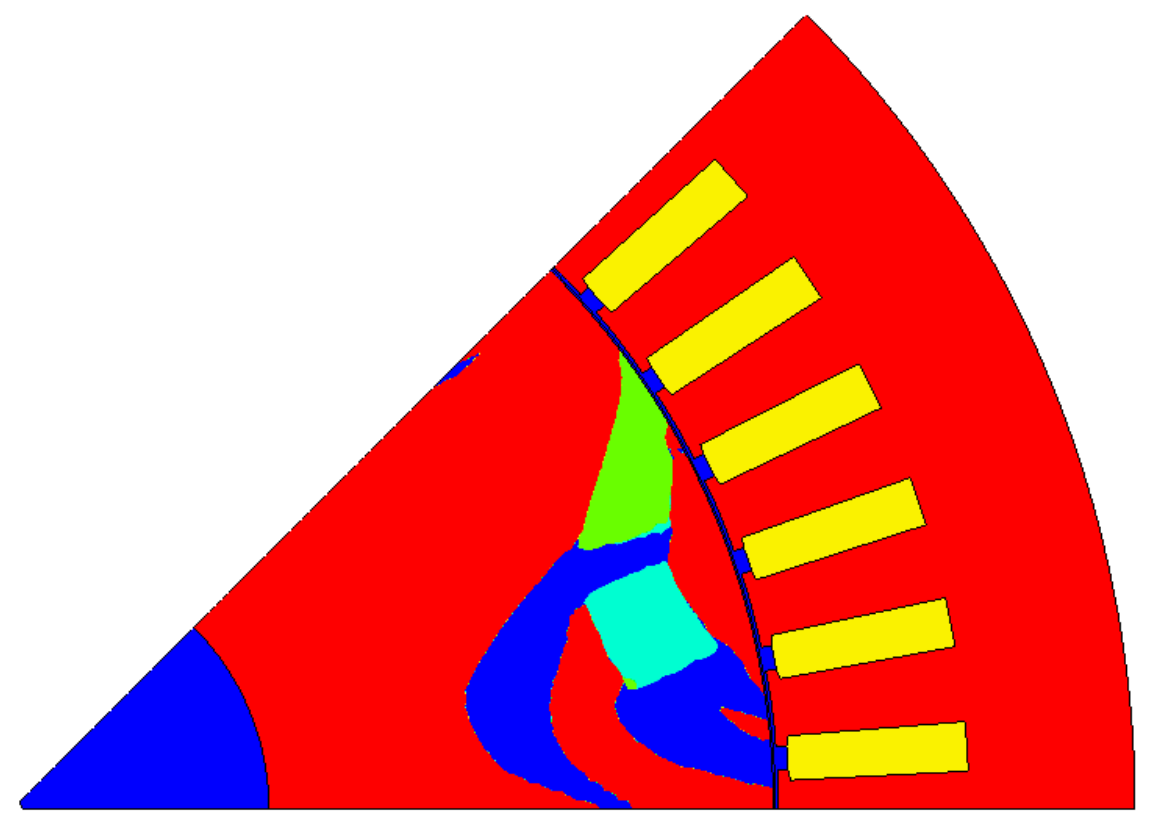} \\
		$\overline{\mathcal{T}(\Omega)}=913\text{Nm}$&$\overline{\mathcal{T}(\Omega)}=922\text{Nm}$&$\overline{\mathcal{T}(\Omega)}=893\text{Nm}$ \\
		\includegraphics[trim=19.7cm 7.3cm 19.7cm 7.3cm,clip,width=0.32\textwidth]{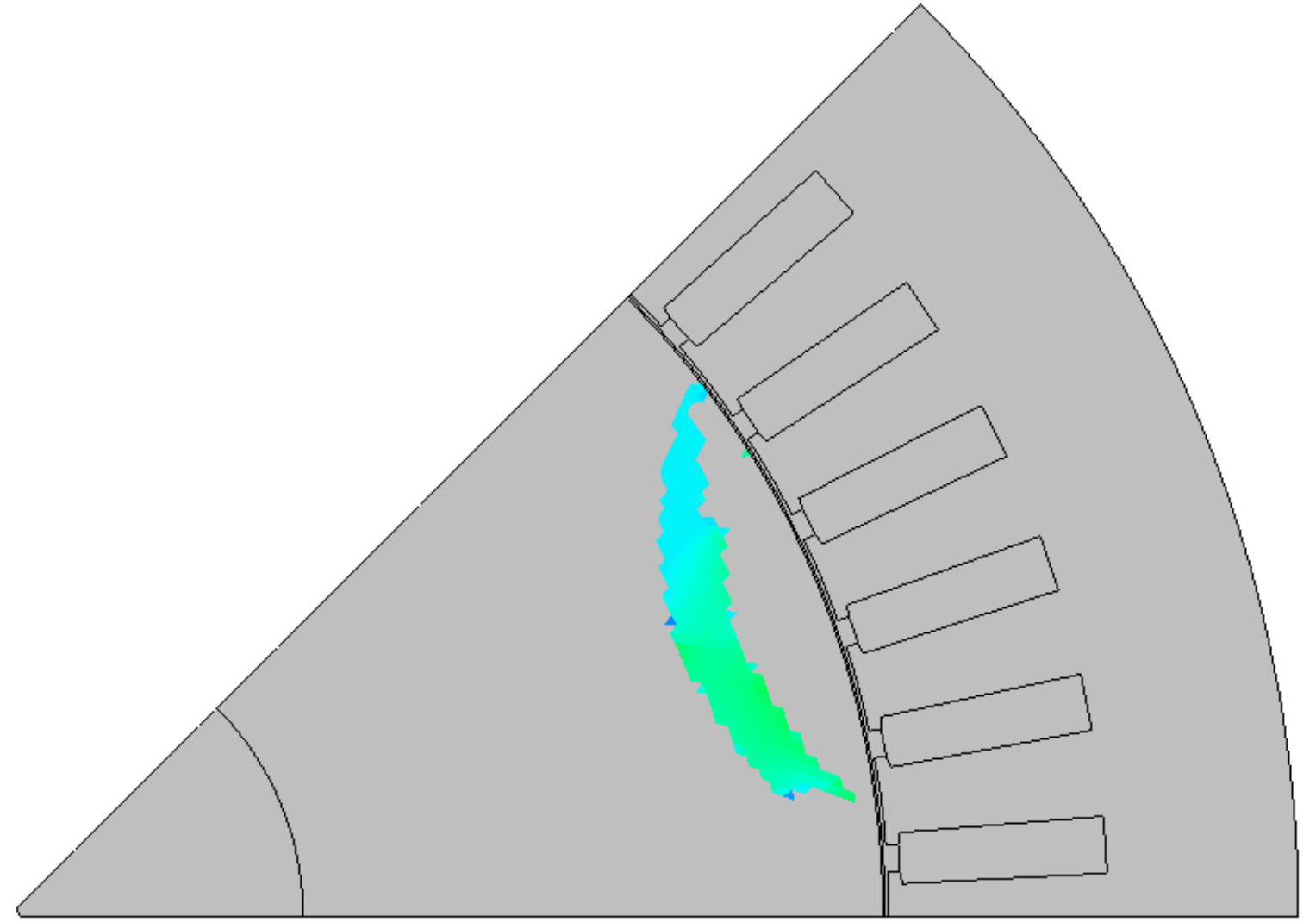} 
		&\includegraphics[trim=19.7cm 7.3cm 19.7cm 7.3cm,clip,width=0.32\textwidth]{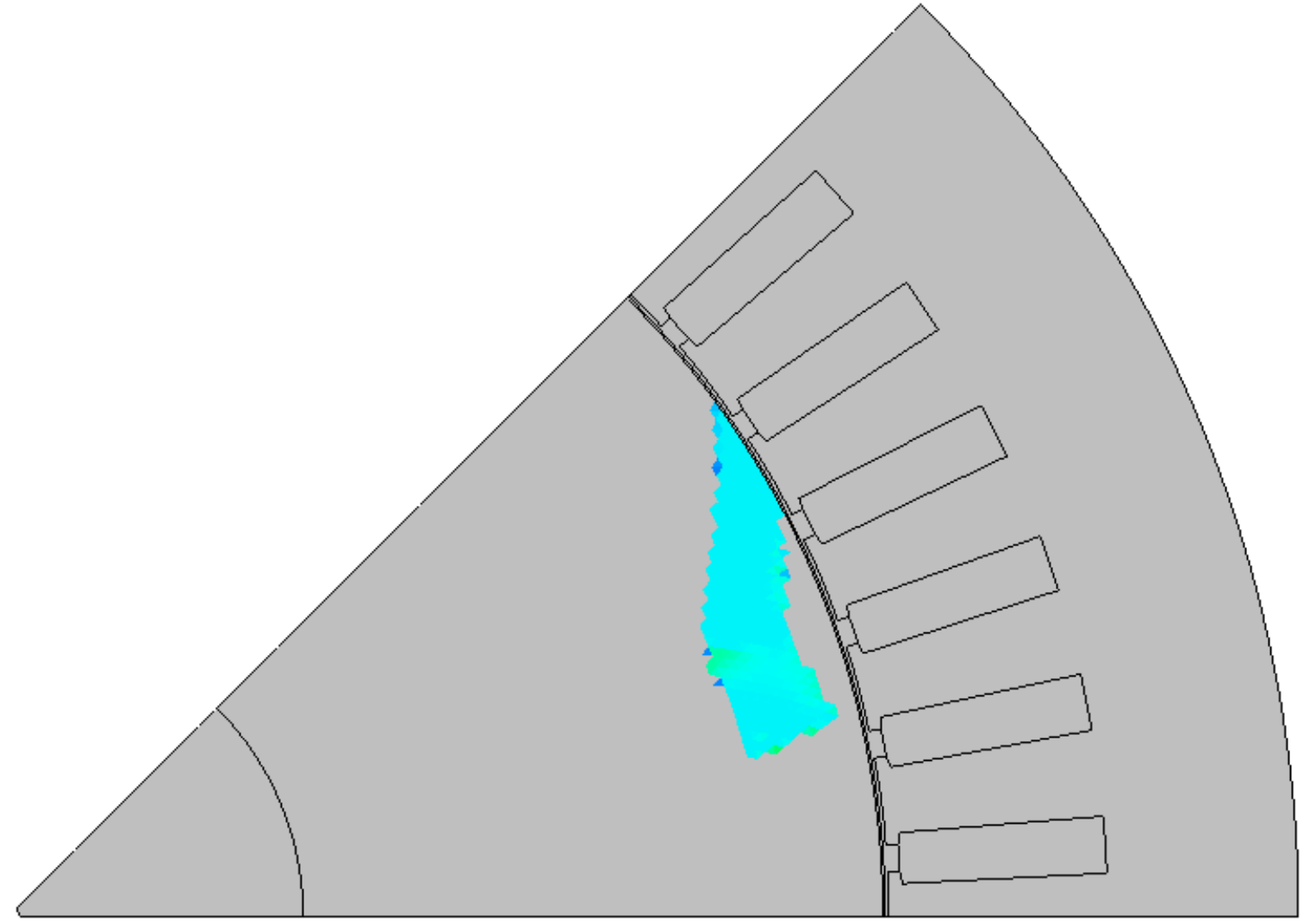} 
		&\includegraphics[trim=19.7cm 7.3cm 19.7cm 7.3cm,clip,width=0.32\textwidth]{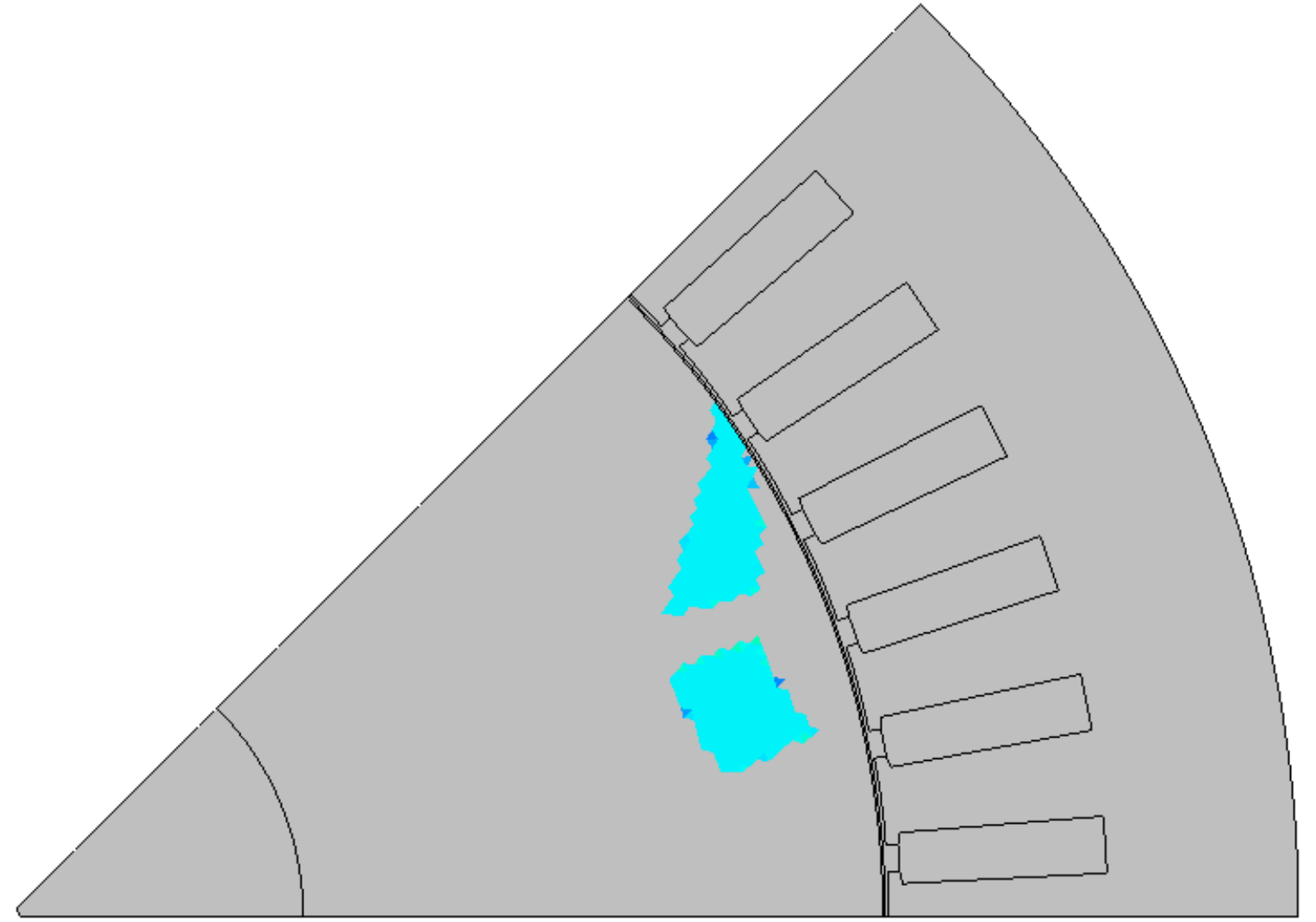} \\
		$\overline{\mathcal{D}(\Omega)}=10\%$&$\overline{\mathcal{D}(\Omega)}=3\%$&$\overline{\mathcal{D}(\Omega)}=3\%$\\
		\includegraphics[trim=19.7cm 7.3cm 19.7cm 7.3cm,clip,width=0.32\textwidth]{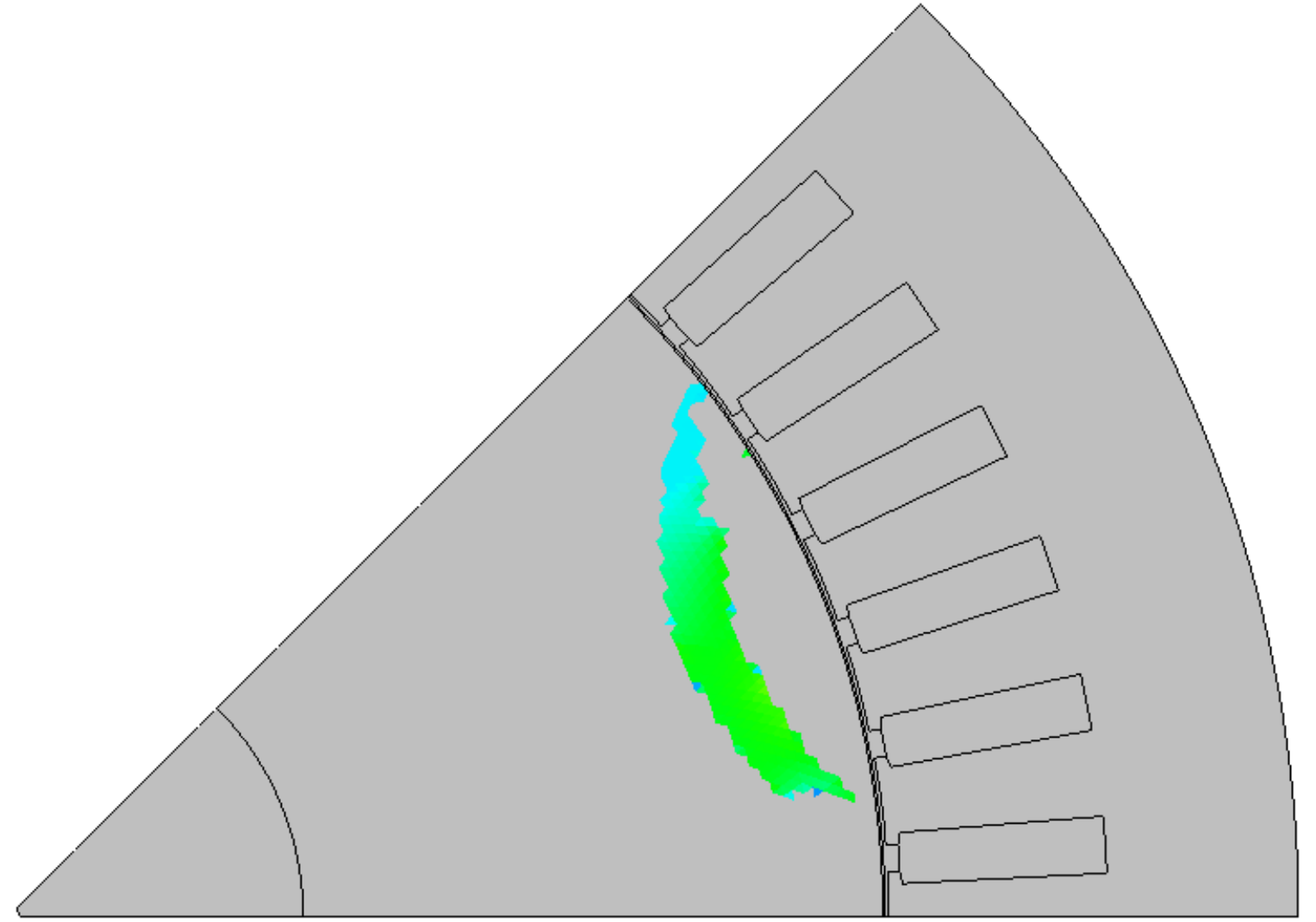} 
		&\includegraphics[trim=19.7cm 7.3cm 19.7cm 7.3cm,clip,width=0.32\textwidth]{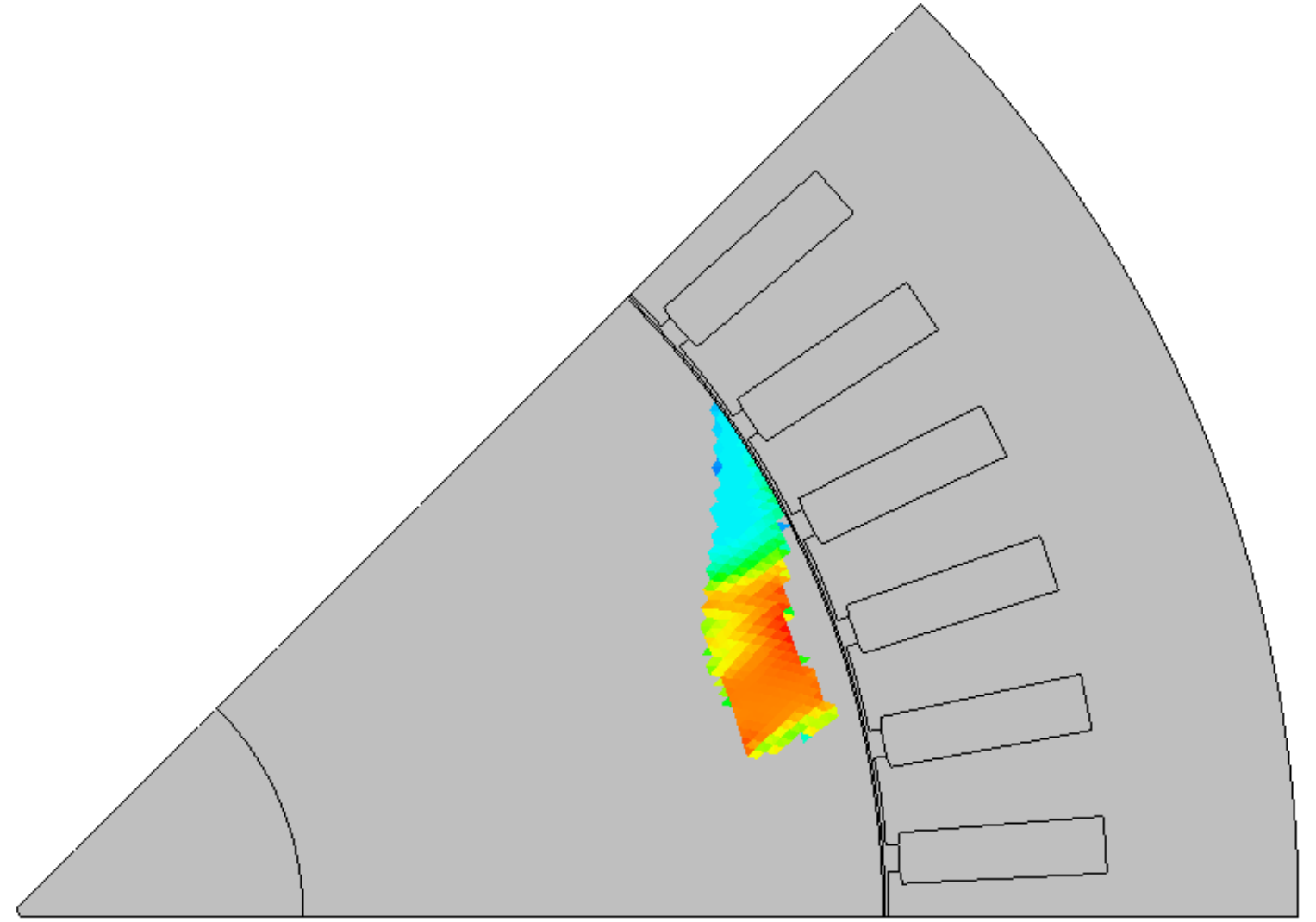} 
		&\includegraphics[trim=19.7cm 7.3cm 19.7cm 7.3cm,clip,width=0.32\textwidth]{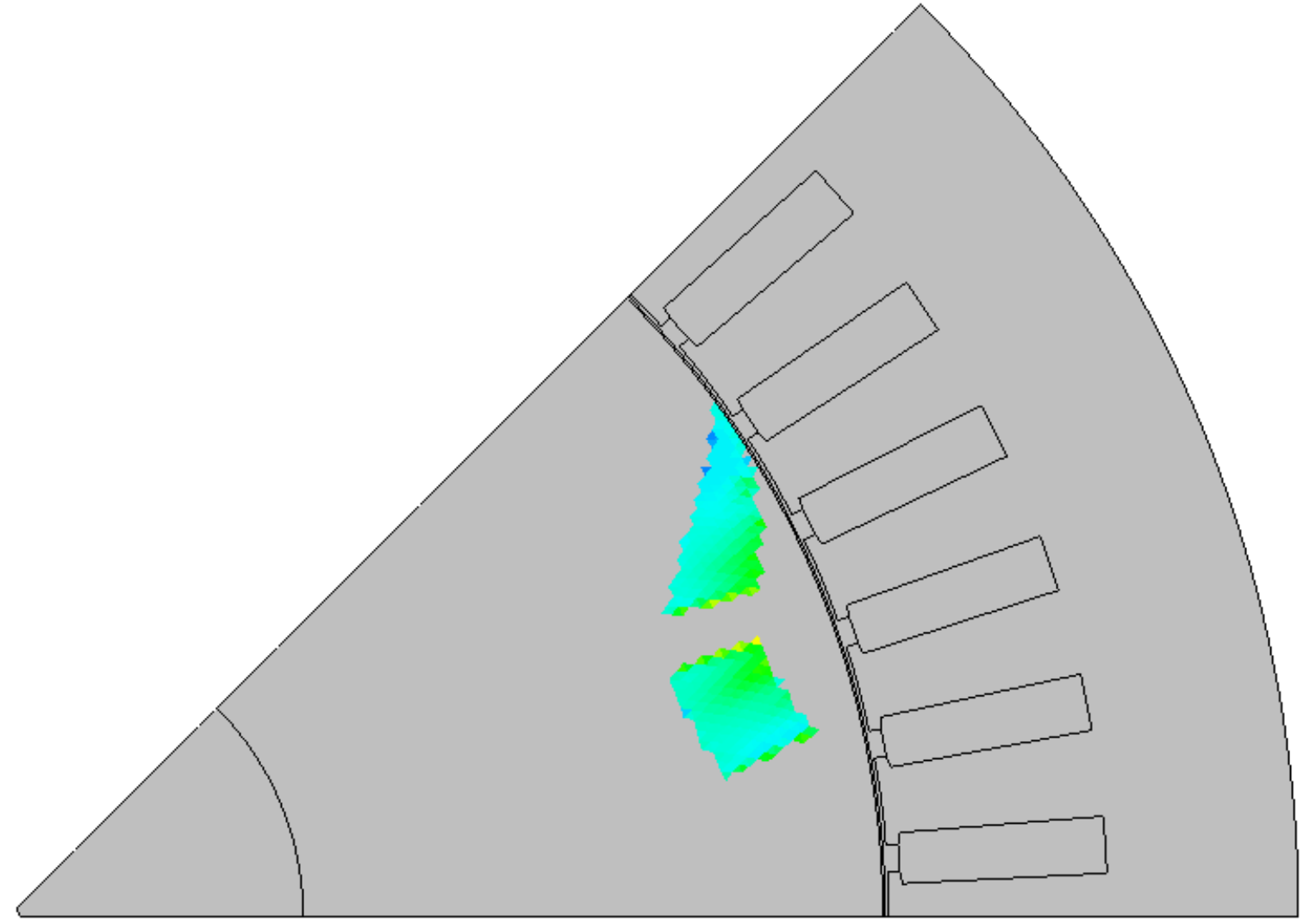} \\
		$\overline{\mathcal{D}(\Omega)}=24\%$&$\overline{\mathcal{D}(\Omega)}=47\%$&$\overline{\mathcal{D}(\Omega)}=13\%$
	\end{tabular}
	\caption{Optimal designs (first line). Partial demagnetization at nominal operation (second line) and at damaging operation (third line), light blue means no demagnetization, green and yellow partial and red full demagnetization. The overline denotes the average for $N=11$ rotor positions of torque and demagnetization resp.}
\end{figure}
\subsection{Discussion}
In all optimizations we maximized the torque while constraining the magnet volume with 10\% of the rotor volume. In the first column we used the linear magnet material law (\ref{krenn:law_linear}) whereas we considered the non-linear law (\ref{krenn:law_nonlinear}) in column two and three. In the third column we added the constraint on demagnetizing fields at the damaging operation point (\ref{krenn:PF_constr})-(\ref{krenn:PF_magnetoD}). We observe in the second line an improvement of the demagnetization robustness from 10\% to 3\% by considering the non-linear magnet which cannot be further improved by the additional constraint. This fits the observations of the authors in \cite{theodor}. If we consider the damaging operation point in the third line the second design demagnetizes heavily. The additional constraint is able to reduce this from 47\% to 13\% while loosing only 3\% of the torque.
\begin{remark}
	The remaining demagnetization is probably due to the fact that we considered a relaxation (\ref{krenn:weak_constr}) of the original constraint. One needs to investigate the role of the parameters $p$ in (\ref{krenn:eq_rho}), $B^*$ and $\gamma$. The problem turned out to be very sensitive to large values of $\gamma$, the results were obtained with $\gamma=10$.
\end{remark}
\section{Outlook}
The next possible step is to model the demagnetization process by adapting the material law after each pseudo-time step and constraining the performance of the demagnetized machine in a multi-stage simulation. Since demagnetization happens also with fields perpendicular to the magnet we aim to include this too by an improved material model.

\section*{Acknowledgments} The authors acknowledge the support via the FWF-funded collaborative research center F90 project 08 (CREATOR – Computational Electric Machine Laboratory). Part of this work was accomplished during an ”Industry Sabbatical” research stay of N. Krenn at Robert Bosch GmbH, Corporate Research in Renningen.

\bibliographystyle{abbrv}

\end{document}